\documentclass[twoside,11pt]{article}
\usepackage{amsfonts}
\usepackage{graphicx,color}
\usepackage{amssymb}
\usepackage{amsmath}

\input{tcilatex}

\begin{document}

\date{}

\centerline{}

\centerline {\Large{\bf Notes on the q-Analogues of the Natural Transforms
and  }}
\centerline {}
\centerline{\Large{\bf  Some  Further Applications }}

\centerline{}

\centerline{\bf {S. K. Q. Al-Omari} }

\centerline{}

\centerline{Department of Applied Sciences, Faculty of Engineering
Technology; Al-Balqa' } \centerline{ Applied University, Amman 11134, Jordan}
\centerline{s.k.q.alomari@fet.edu.jo}

\centerline{}

\centerline{\bf {Adem K{\i}l{\i}\c{c}man} }

\centerline{}

\centerline{Department of Mathematics, University Putra Malaysia (UPM)}
\centerline{43400 Serdang, Selangor, Malaysia}
\centerline{akilic@upm.edu.my }

\centerline{}\quad

\centerline{\bf Abstract}

\centerline{\bf }

As an extension to the Laplace and Sumudu transforms the classical Natural
transform was proposed to solve certain fluid flow problems. In this paper,
we investigate $q$-analogues of the $q$-Natural transform of some special
functions. We derive the $q$-analogues of the $q$-integral transform and
further apply to some general special functions such as : the exponential
functions, the $q$-trigonometric functions, the $q$-hyperbolic functions and
the Heaviside Function. Some further results involving convolutions and
differentiations are also obtained.

\begin{description}
\item[\textbf{Keywords:}] $q$-hyperbolic function; $q$-trigonometric
function, $q$-Natural transform; Heaviside Function; Natural transform.
\end{description}


\section{Introduction}

The subject of fractional calculus (integrals and derivatives of any real or
complex order) has gained noticeable importance and popularity due to mainly
its demonstrated applications in many seemingly diverse fields of science
and engineering. Much of the theory of fractional calculus is based upon the
familiar Riemann-Liouville fractional derivatives and integrals. Recently,
there was a significant increase of activity in the area of the q-calculus
due to applications of the $q$-calculus in mathematics, statistics and
physics.\\

Jackson in \cite{11} presented a precise definition of the so-called $q$%
-Jackson integral and developed a $q$-calculus in a systematic way. Some
remarkable integral transforms have different $q$-analogues in the theory of
$q$-calculus. Among those $q$-integrals we recall here are the $q$-Laplace
integral transform \cite{2,3,4}, $q$-Sumudu integral transform \cite{5,6}, Weyl fractional $q$-integral operator \cite{7},
$q$-Wavelet integral transform \cite{8}, $q$-Mellin integral
transform \cite{9}, and few others. In this paper, we give some $q$%
-analogues of some recently investigated transform named as Natural
transform and obtain some desired $q$-properties.\\

In the following section, we present some notations and terminologies from
the $q$-calculus. In Section 3, we recall the definition and properties of
the Natural transform. In Section 4, we derive the definition of the first $%
q $-analogue of the Natural transform and apply to some special functions.
Sections 5-7 are devoted to some applications of the first $q$-analogue of
the $q$-Natural transform to Heaviside Functions, convolutions and
differentiations. The remaining two sections are investigating the second $q$%
-analogue of the $q$-Natural transform of some elementary functions and some
applications.

\section{Definitions and Preliminaries}

\noindent We present some usual notions and notations used in the $q$%
-calculus see \cite{10,11,12}. Throughout
this paper, we assume $q$ to be a fixed number satisfying $0<q<1$.\\

\noindent The $q$-calculus beings with the definition of the $q$-analogue $%
d_{q}f\left( x\right) $ of the differential of functions,%
\begin{equation}
d_{q}f\left( x\right) =f\left( qx\right) -f\left( x\right) .  \tag{1}
\end{equation}%
Having said this, we immediately get the $q$-analogue of the derivative of $%
f\left( x\right) $, called its $q$-derivative,%
\begin{equation}
\left( D_{q}f\right) \left( x\right) :=\frac{d_{q}f\left( x\right) }{d_{q}x}%
:=\frac{f\left( x\right) -f\left( qx\right) }{\left( 1-q\right) x}, \ \text{ if }%
x\neq 0,  \tag{2}
\end{equation}%
$\left( D_{q}f\right) \left( 0\right) =\acute{f}\left( 0\right) $, provided $%
\acute{f}\left( 0\right) $ exists. If $f$ is differentiable, then $\left(
D_{q}f\right) \left( x\right) $ tends to $\acute{f}\left( 0\right) $ as $q$
tends to $1$.\\

\noindent Notice that the $q$-derivative satisfies the following $q$%
-analogue of Leibnitz rule,%
\begin{equation}
D_{q}\left( f\left( x\right) g\left( x\right) \right) =g\left( x\right)
D_{q}f\left( x\right) +f\left( qx\right) D_{q}g\left( x\right) .  \tag{3}
\end{equation}%
The $q$-Jackson integrals from $0$ to $x$ and from $x$ to $\infty $ are
defined in \cite{9, 11} by
\begin{equation}
\dint\nolimits_{0}^{x}f\left( x\right) d_{q}x=\left( 1-q\right)
x\sum_{0}^{\infty }f\left( xq^{k}\right) q^{k},  \tag{4}
\end{equation}

\begin{equation}
\dint\nolimits_{0}^{\infty }f\left( x\right) d_{q}x=\left( 1-q\right)
x\sum_{-\infty }^{\infty }f\left( q^{k}\right) q^{k},  \tag{5}
\end{equation}%
\noindent provided the sum converges absolutely.\\

\noindent The $q$-Jackson integral in a generic interval $\left[ a,b\right] $
is given by in \cite{11} as
\begin{equation}
\dint\nolimits_{a}^{b}f\left( x\right) d_{q}x=\dint\nolimits_{0}^{b}f\left(
x\right) d_{q}x-\dint\nolimits_{0}^{a}f\left( x\right) d_{q}x.  \tag{6}
\end{equation}%
The improper integral is defined in the way that%
\begin{equation}
\dint\nolimits_{0}^{\frac{\infty }{A}}f\left( x\right) d_{q}x=\left(
1-q\right) \sum_{-\infty }^{\infty }f\left( \frac{q^{k}}{A}\right) \frac{%
q^{k}}{A}  \tag{7}
\end{equation}%
and, for $n\in \mathbb{Z}$, we have%
\begin{equation}
\dint\nolimits_{0}^{\frac{\infty }{q^{n}}}f\left( x\right)
d_{q}x=\dint\nolimits_{0}^{\infty }f\left( x\right) d_{q}x.  \tag{8}
\end{equation}%
The $q$-integration by parts is defined for functions $f$ and $g$ by%
\begin{equation}
\dint\nolimits_{0}^{b}g\left( x\right) D_{q}f\left( x\right) d_{q}x=f\left(
b\right) g\left( b\right) -f\left( a\right) g\left( a\right)
-\dint\nolimits_{a}^{b}f\left( qx\right) D_{q}g\left( x\right) d_{q}x.
\tag{9}
\end{equation}%
For $x\in \mathbb{C}$, the $q$-shifted factorials are defined by%
\begin{equation*}
\left( x;q\right) _{0}=1;\left( x,q\right) _{t}=\frac{\left( x;q\right)
_{\infty }}{\left( aq^{x};q\right) _{\infty }};\left( x,q\right)
_{n}=\prod\limits_{k=0}^{n-1}\left( 1-xq^{k}\right)
\end{equation*}%
and%
\begin{equation}
\left( x;q\right) _{\infty }=\prod\limits_{k=0}^{\infty }\left(
1-xq^{k}\right) ,  \tag{10}
\end{equation}%
$n=0,1,2,...$.\\

\noindent The $q$-analogue of $x$ and $\infty $ is defined by%
\begin{equation}
\left[ x\right] =\frac{1-q^{x}}{1-q}\text{ and }\left[ \infty \right] =\frac{%
1}{1-q}.  \tag{11}
\end{equation}%
\noindent The important $q$-analogues of the exponential function of first
and second kinds are respectively given as:%
\begin{equation}
E_{q}\left( x\right) =\sum_{n=0}^{\infty }q^{\frac{n\left( n-1\right) }{2}}%
\frac{x^{n}}{\left( \left[ n\right] _{q}\right) !}\text{\ }\left( t\in
\mathbb{C} \right) ,  \tag{12}
\end{equation}%
and%
\begin{equation}
e_{q}\left( x\right) =\sum_{n=0}^{\infty }\frac{x^{n}}{\left( \left[ n\right]
_{q}\right) !}\text{\ }\left( \left\vert t\right\vert <1\right) ,  \tag{13}
\end{equation}%
where $\left( \left[ n\right] _{q}\right) !=\left[ n\right] _{q}\left[ n-1%
\right] _{q}...\left[ 2\right] _{q}\left[ 1\right] _{q},\left[ n\right] =%
\dfrac{1-q^{n}}{1-q}=q^{n-1}+...+q+1$.\\

\noindent However, due to the product expansions, $\left( e_{q}\left(
x\right) \right) ^{-1}=E_{q}\left( -x\right) \left( \text{not }e_{q}\left(
-x\right) \right) $, which explains the need of both $q$-analogues of the
exponential function.\\

\noindent The $q$-derivative of $E_{q}\left( x\right) $ is $D_{q}E_{q}\left(
x\right) =E_{q}\left( qx\right) $, whereas, the $q$-derivative of $%
e_{q}\left( x\right) $ is $D_{q}e_{q}\left( x\right) =e_{q}\left( x\right)
,e_{q}\left( 0\right) =1$.\\

\noindent The gamma and beta functions satisfy the $q$-integral
representations%
\begin{equation}
\left.
\begin{array}{l}
\Gamma _{q}\left( t\right) =\dint\nolimits_{0}^{\frac{1}{1-q}%
}x^{t-1}E_{q}\left( -qx\right) d_{q}x \\
\text{and} \\
B_{q}\left( t;s\right) =\dint\nolimits_{0}^{1}x^{t-1}\left( 1-qx\right)
_{q}^{s-1}d_{q}x,\left( t,s>0\right)%
\end{array}%
\right\}  \tag{14}
\end{equation}%
that satisfy $B_{q}\left( t;s\right) =\dfrac{\Gamma _{q}\left( s\right)
\Gamma _{q}\left( t\right) }{\Gamma _{q}\left( s+t\right) }$ and $\Gamma
_{q}\left( t+1\right) =\left[ t\right] _{q}\Gamma _{q}\left( t\right) $.
Due to $\left( 12\right) $ and $\left( 13\right) $, the $q$-analogues of
sine and cosine functions of the second and first kinds are respectively
given as :
\begin{equation}
\left.
\begin{array}{l}
\sin ^{q}\left( at\right) =\dsum_{0}^{\infty }\left( -1\right) ^{n}\frac{%
\left( at\right) ^{2n+1}}{\left( \left[ 2n+1\right] _{q}\right) !}; \\
\cos ^{q}\left( at\right) =\dsum_{0}^{\infty }\left( -1\right) ^{n}\frac{%
\left( at\right) ^{2n}}{\left( \left[ 2n\right] _{q}\right) !}; \\
\sin _{q}\left( at\right) =\dsum_{0}^{\infty }\left( -1\right) ^{n}\frac{q^{%
\frac{n\left( n+1\right) }{2}}}{\left( \left[ 2n+1\right] _{q}\right) !}%
a^{2n+1}t^{2n+1}; \\
\cos _{q}\left( at\right) =\dsum_{0}^{\infty }\left( -1\right) ^{n}\frac{q^{%
\frac{n\left( n-1\right) }{2}}}{\left( \left[ 2n\right] _{q}\right) !}%
a^{2n}t^{2n}.%
\end{array}%
\right\}   \tag{15}
\end{equation}

\section{The Natural Transform}

\noindent The Natural transform of a function $f\left( x\right)$ on $0<x<\infty$ then it was proposed by Khan and Khan \cite{13} as an extension to the
Laplace and Sumudu transforms to solve some fluid flow problems.\\

\noindent Later, Silambarasan and Belgacem \cite{16} have derived
certain electric field solutions of the Maxwell's equation in conducting
media. In \cite{17}, the author applied the Natural transform to
some ordinary differential equations and some space of Boehmians. Further
investigation of the Natural transform can be obtained from \cite{14} and \cite{15}.\\

\noindent The Natural transform of a function $f\left( t\right), 0<t<\infty$ is defined over the set $A$, where%
\begin{equation*}
A=\left\{ f\left( t\right) \left\vert \exists M,\tau _{1},\tau
_{2}>0,\left\vert f\left( t\right) \right\vert <Me^{t/\tau j},\text{ if }%
t\in \left( -1\right) ^{j}\times \left[ 0,\infty \right) \right\vert \right\}
\end{equation*}%
by $\left( \cite{13},\left( 1\right) \right) $%
\begin{equation}
\left( Nf\right) \left( u;v\right) =\dint\nolimits_{0}^{\infty }f\left(
ut\right) \exp \left( -vt\right) dt\text{ \ }\left( u,v>0\right) .  \tag{16}
\end{equation}%
Provided the integral on the right of $\left( 16\right) $ exists it is easy
to see that%
\begin{equation}
\left( Nf\right) \left( u;1\right) =\left( Sf\right) \left( u\right) \text{
and }\left( Nf\right) \left( 1,v\right) =\left( Lf\right) \left( v\right)
\tag{17}
\end{equation}%
where $Sf$ and $Lf$ are respectively the Sumudu and Laplace transforms of $%
f$.\\

\noindent Moreover, the Natural-Laplace and Natural-Sumudu dualities are
given in \cite{17,5,6} as
\begin{equation}
\left( Nf\right) \left( u;v\right) =\frac{1}{u}\dint\nolimits_{0}^{\infty
}f\left( t\right) \exp \left( \frac{-vt}{u}\right) dt  \tag{18}
\end{equation}%
and%
\begin{equation}
\left( Nf\right) \left( u;v\right) =\frac{1}{v}\dint\nolimits_{0}^{\infty
}f\left( \frac{ut}{v}\right) \exp \left( -t\right) dt,  \tag{19}
\end{equation}%
respectively.\\

\noindent It further from $\left( 18\right) $ and $\left( 19\right) $ can be
easily observed that%
\begin{equation}
\left( Nf\right) \left( u;v\right) =\frac{1}{u}\left( Lf\right) \left( \frac{%
u}{v}\right) \text{ and }\left( Nf\right) \left( u;v\right) =\frac{1}{v}%
\left( Sf\right) \left( \frac{u}{v}\right) .  \tag{20}
\end{equation}%
Some values of the Natural transform of some known functions we mention here
are $\left[ \cite{17},\text{p}.731\right] $

$\left( i\right) N\left( a\right) \left( u;v\right) =\dfrac{1}{v}$, where $a$
is a constant.

$\left( ii\right) N\left( \delta \right) \left( u;v\right) =\dfrac{1}{v}$,
where $\delta $ is the delta function.

$\left( iii\right) N\left( e^{at}\right) \left( u;v\right) =\dfrac{1}{v-au}, a $ is a constant.

$\left( iv\right) $ The scaling property is written in two ways as%
\begin{equation*}
\left( Nf\left( kt\right) \right) \left( u;v\right) =\frac{1}{k}\left(
Nf\right) \left( ku;v\right) \text{ and }\left( Nf\left( kt\right) \right)
\left( u;v\right) =\frac{1}{k}\left( Nf\right) \left( u;\frac{v}{k}\right) .
\end{equation*}

\section{The $q$-Analogue of the $q$-Natural Transform of First Kind}

\noindent Hahn \cite{18} and later Ucar and Albayrak \cite{4} defined the $q$-analogue of first and second types of the
well-known Laplace transform by means of the $q$-integrals%
\begin{equation}
L_{q}\left( f\left( t\right) ;s\right) =\frac{1}{1-q}\dint\nolimits_{0}^{%
\frac{1}{s}}E_{q}\left( qst\right) f\left( t\right) d_{q}t  \tag{21}
\end{equation}%
and%
\begin{equation}
_{q}L\left( f\left( t\right) ;s\right) =\frac{1}{1-q}\dint\nolimits_{0}^{%
\infty }e_{q}\left( -st\right) f\left( t\right) d_{q}t.  \tag{22}
\end{equation}%
The $q$-analogues of the Sumudu transform of first and second types are
defined by \cite{5,6}
\begin{equation}
S_{q}\left( f\left( t\right) ;s\right) =\frac{1}{\left( 1-q\right) _{s}}%
\dint\nolimits_{0}^{s}E_{q}\left( \frac{q}{s}t\right) f\left( t\right) d_{q}t
\tag{23}
\end{equation}%
and%
\begin{equation}
_{q}S\left( f\left( t\right) ;s\right) =\frac{1}{1-q}\dint\nolimits_{0}^{%
\infty }e_{q}\left( -\frac{t}{s}\right) f\left( t\right) d_{q}t  \tag{24}
\end{equation}%
where $s\in \left( -\tau _{1},\tau _{2}\right) $ and$\ f$ is a function
belongs to the set $A$,%
\begin{equation*}
A=\left\{ f\left( t\right) \left\vert \exists M,\tau _{1},\tau
_{2}>0,\left\vert f\left( t\right) \right\vert <Me^{\frac{\left\vert
t\right\vert }{\tau _{j}}},t\in \left( -1\right) ^{j}\times \left[ 0,\infty
\right) \right. \right\} .
\end{equation*}%
Now, we are in a position to demonstrate our results as follows. The $q$%
-analogue of the Natural transform of first kind as%
\begin{equation}
\left( N_{q}f\right) \left( u;v\right) :=\left( N_{q}f;\left( u,v\right)
\right) :=\frac{1}{u}\dint\nolimits_{0}^{\infty }E_{q}\left( -q\frac{vt}{u}%
\right) d_{q}t,  \tag{25}
\end{equation}%
provided the function $f\left( t\right) $ is defined on $A$ and, $u$ and $v$
are the transform variables.

\noindent The series representation of $\left( 25\right) $ can be written as:%
\begin{equation}
\left( N_{q}f\right) \left( u;v\right) =\frac{1}{\left( 1-q\right) u}%
\sum_{k\in
\mathbb{Z}
}q^{k}f\left( q^{k}\right) E_{q}\left( -q^{k+1}\frac{v}{u}\right),  \tag{26}
\end{equation}%
and by $\left( 10\right) $, $\left( 26\right) $ can be put into the form%
\begin{equation}
\left( N_{q}f\right) \left( u;v\right) =\frac{\left( q;q\right) _{\infty }}{%
\left( 1-q\right) ^{n}}\sum_{k\in
\mathbb{Z}
}\frac{q^{k}f\left( q^{k}\right) }{\left( -\dfrac{v}{u};q\right) _{k+1}}.
\tag{27}
\end{equation}%
We derive now some values of $N_{q}$ of some special functions. \\

\noindent \textbf{Theorem 1}. Let $\alpha \in \mathbb{R}$, then we have%
\begin{equation}
\left( N_{q}t^{\alpha }\right) \left( u;v\right) =\frac{u^{\alpha }}{%
v^{\alpha +1}}\Gamma \left( \alpha +1\right) .  \tag{28}
\end{equation}%
\noindent \textbf{Proof}. By setting the variables, we get%
\begin{eqnarray*}
\left( N_{q}t^{\alpha }\right) \left( u;v\right)  &=&\frac{u^{\alpha }}{%
v^{\alpha +1}}\dint\nolimits_{0}^{\infty }t^{\alpha }E_{q}\left( -qt\right)
d_{q}t \\
&=&\frac{u^{\alpha }}{v^{\alpha +1}}\Gamma _{q}\left( \alpha +1\right) .
\end{eqnarray*}%
The theorem hence follows. \\

\noindent A direct corollary of $\left( 28\right) $ can be%
\begin{equation}
\left( N_{q}t^{n}\right) \left( u;v\right) =\frac{u^{n}}{v^{n+1}}\left( %
\left[ n\right] _{q}\right) !.  \tag{29}
\end{equation}

\noindent \textbf{Lemma 2}. Let $a$ be a positive real number. Then, we have%
\begin{equation}
\left( N_{q}e_{q}\left( at\right) \right) \left( u;v\right) =\dfrac{1}{v-au},%
\text{ }au<v.  \tag{30}
\end{equation}%
\noindent \textbf{Proof}. Using the second kind $q$-analogue of the $q$%
-exponential function, we write%
\begin{equation}
\left( N_{q}e_{q}\left( at\right) \right) \left( u;v\right) =\frac{1}{u}%
\dint\nolimits_{0}^{\infty }e_{q}\left( at\right) E_{q}\left( -q\frac{v}{u}%
t\right) d_{q}t.  \tag{31}
\end{equation}%
In series representation $\left( 31\right) $ can be written as%
\begin{equation}
\left( N_{q}e_{q}\left( at\right) \right) \left( u;v\right)
=\sum_{0}^{\infty }\frac{a^{n}}{u\left( \left[ n\right] _{q}\right) !}%
t^{n}E^{q}\left( -q\frac{v}{u}t\right) d_{q}t  \tag{32}
\end{equation}%
By $\left( 29\right) ,\left( 32\right) $ gives a geometric series expansion
and hence,%
\begin{equation*}
\left( N_{q}e_{q}\left( at\right) \right) \left( u;v\right) =\frac{1}{v}%
\sum_{0}^{\infty }\left( \frac{au}{v}\right) ^{n}=\frac{1}{v-au},\text{ }%
au<v.
\end{equation*}%
This completes the proof of the lemma. \\

\noindent \textbf{Theorem 3}. Let $a$ be a positive real number. Then, we
have%
\begin{equation}
\left( N_{q}E_{q}\left( at\right) \right) \left( u;v\right) =\frac{1}{v}%
\sum_{0}^{\infty }q\frac{n\left( n-1\right) }{2}\left( \frac{au}{v}\right)
^{n}.  \tag{33}
\end{equation}

\noindent \textbf{Proof}. By the first kind $q$-analogue of the exponential
function we indeed get%
\begin{equation}
\left( N_{q}E_{q}\left( at\right) \right) \left( u;v\right)
=\sum_{0}^{\infty }\frac{a_{n}q\frac{n\left( n-1\right) }{2}}{\left( \left[ n%
\right] _{q}\right) !}\frac{1}{u}\dint\nolimits_{0}^{\infty
}t^{n}E_{q}\left( -q\frac{v}{u}t\right) d_{q}t.  \tag{34}
\end{equation}

\noindent The parity of $\left( 29\right) $ gives%
\begin{equation*}
\left( N_{q}E_{q}\left( at\right) \right) \left( u;v\right) =\frac{1}{v}\sum
q\frac{n\left( n-1\right) }{2}\left( \frac{au}{v}\right) ^{n}.
\end{equation*}

\noindent This completes the proof of the theorem. \\

\noindent The hyperbolic $q$-cosine and $q$-sine functions are given as $$\cosh ^{q}t=\dfrac{e_{q}\left( t\right) +e_{q}\left( -t\right) }{2} \ {\rm \ and} \ \
\sinh ^{q}t=\dfrac{e_{q}\left( t\right) -e_{q}\left( -t\right) }{2}.$$ Hence,
as a corollary of Theorem 2, we have%
\begin{eqnarray}
\left( N_{q}\cosh ^{q}t\right) \left( u;v\right)  &=&\frac{1}{2}\left\{
\left( N_{q}e_{q}\left( at\right) \right) \left( u,v\right) +\left(
N_{q}e_{q}\left( -t\right) \right) \left( u,v\right) \right\}   \notag \\
&=&\frac{v}{v^{2}+a^{2}u^{2}},\text{ }au<v.  \TCItag{35}
\end{eqnarray}%
and%
\begin{equation}
\left( N_{q}\sinh ^{q}t\right) \left( u;v\right) =\frac{au}{v^{2}+a^{2}u^{2}}%
,\text{ }au<v.  \tag{36}
\end{equation}%
\noindent \textbf{Theorem 4}. Let $a$ be a positive real number. Then, we
have%
\begin{equation*}
\left( N_{q}\cos ^{q}at\right) \left( u;v\right) =\frac{v}{v^{2}-a^{2}u^{2}},
\end{equation*}%
provided $au<v$.

\noindent \textbf{Proof}. On account of $\left( 15\right) $ and $\left(
29\right) $ we obtain%
\begin{eqnarray*}
\left( N_{q}\cos ^{q}at\right) \left( u;v\right)  &=&\sum_{0}^{\infty }\frac{%
a^{2n}}{u\left( \left[ 2n\right] _{q}\right) !}\dint\nolimits_{0}^{\infty
}t^{2n}E_{q}\left( -q\frac{v}{u}t\right) d_{q}t \\
&=&\frac{1}{v}\sum_{0}^{\infty }\left( \frac{au}{v}\right) ^{2n}.
\end{eqnarray*}%
For $au<v$, the geometric series converges to the sum%
\begin{equation*}
\left( N_{q}\cos ^{q}at\right) \left( u;v\right) =\frac{v}{v^{2}-a^{2}u^{2}}%
\text{ provided }au<v.
\end{equation*}%
\noindent This completes the proof of the theorem. \\

\noindent Similarly, by $\left( 15\right) $ and $\left( 29\right) $,\ we
deduce that%
\begin{equation}
\left( N_{q}\sin ^{q}at\right) \left( u;v\right) =\frac{au}{v^{2}-a^{2}u^{2}}%
,\text{ }au<v.  \tag{38}
\end{equation}

\section{$N_{q}$ and $q$-Differentiation}

\noindent In this section of this paper we discuss some $q$-differentiation
formulae.

\noindent On account of $\left( 12\right) $ we derive the following
differentiation result.

\noindent \textbf{Lemma 5}. Let $u,v>0$, then we have%
\begin{equation}
D_{q}E_{q}\left( -q\frac{u}{v}t\right) =\frac{v}{u}\sum_{0}^{\infty }\left(
-1\right) ^{n+1}q^{\frac{\left( n+1\right) \left( n+2\right) }{2}}\frac{v^{n}%
}{u^{n}}t^{n}.  \tag{39}
\end{equation}

\noindent \textbf{Proof}. By using the $q$-representation of $E_{q}$ in $%
\left( 12\right) $ we write%
\begin{eqnarray*}
D_{q}E_{q}\left( -q\frac{v}{u}t\right)  &=&D_{q}\sum_{0}^{\infty }\frac{%
\left( -1\right) ^{\frac{n\left( n-1\right) }{2}}}{\left( \left[ n\right]
_{q}\right) !}\left( \frac{qv}{u}\right) ^{n}t^{n} \\
&=&\sum_{1}^{\infty }\left( -1\right) ^{n}\frac{q^{\frac{\left( n+1\right)
\left( n+2\right) }{2}}}{\left( \left[ n-1\right] _{q}\right) !}q^{n}\frac{%
v^{n}}{u^{n}}t^{n-1} \\
&=&\sum_{1}^{\infty }\left( -1\right) ^{n+1}q^{\frac{\left( n+1\right)
\left( n+2\right) }{2}}\frac{v^{n+1}}{u^{n+1}}t^{n}.
\end{eqnarray*}%
Hence, it follows that%
\begin{equation*}
D_{q}E_{q}\left( -q\frac{v}{u}t\right) =\frac{v}{u}\sum_{0}^{\infty }\left(
-1\right) ^{n+1}q^{\frac{\left( n+1\right) \left( n+2\right) }{2}}\frac{v^{n}%
}{u^{n}}t^{n}.
\end{equation*}

\noindent This completes the proof of the theorem.\\

\noindent The Natural transform of the $q$-derivative $D_{q}f$ can be
written as follows.

\noindent \textbf{Theorem 6}. Let $u,v>0$, then we have%
\begin{equation}
N_{q}\left( D_{q}f\left( t\right) \right) \left( u;v\right) =-f\left(
0\right) +\frac{v}{u}N_{q}\left( f\right) \left( u;v\right).  \tag{40}
\end{equation}

\noindent \textbf{Proof}. Using the idea of $q$-integration by parts and the
formula in $\left( 9\right) $ we write%
\begin{eqnarray*}
N_{q}\left( D_{q}f\left( t\right) \right) \left( u;v\right)
&=&\dint\nolimits_{0}^{\infty }D_{q}f\left( t\right) E_{q}\left( -q\frac{v}{u%
}t\right) d_{q}t \\
&=&f\left( t\right) D_{q}E_{q}\left( -q\frac{v}{u}t\right) \QATOPD\vert .
{^{\infty }}{_{0}}-\dint\nolimits_{0}^{\infty }f\left( qt\right)
D_{q}E_{q}\left( -q\frac{v}{u}t\right) d_{q}t
\end{eqnarray*}%
\noindent The parity of Lemma 5 $\left( \text{Eq.}4.39\right) $ gives%
\begin{eqnarray*}
N_{q}\left( D_{q}f\left( t\right) \right) \left( u;v\right) &=&-f\left(
0\right) -\dint\nolimits_{0}^{\infty }f\left( qt\right) \frac{v}{u}%
\sum_{0}^{\infty }\frac{\left( -1\right) ^{n+1}}{\left( \left[ n\right]
_{q}\right) !} \\
&=&-f\left( 0\right) +\frac{v}{u}\dint\nolimits_{0}^{\infty }f\left(
qt\right) \sum_{0}^{\infty }\frac{\left( -1\right) ^{n}}{\left( \left[ n%
\right] _{q}\right) !}q^{\frac{\left( n+1\right) \left( n+2\right) }{2}}%
\frac{v^{n}}{u^{n}}t^{n}d_{q}t
\end{eqnarray*}%
\noindent Changing the variables $qt=y$, and $t^{n}=q^{-n}y^{n}$ imply%
\begin{eqnarray}
N_{q}\left( D_{q}f\left( t\right) \right) \left( u;v\right) &=&-f\left(
0\right) +\frac{v}{u}\dint\nolimits_{0}^{\infty }f\left( y\right)
\sum_{0}^{\infty }\left( -1\right) ^{n}q^{\frac{n\left( n-1\right) }{2}}%
\frac{v^{n}}{u^{n}}y^{n}d_{q}y  \notag \\
&=&-f\left( 0\right) +\frac{v}{u}\dint\nolimits_{0}^{\infty }f\left(
t\right) \sum_{0}^{\infty }\left( -1\right) ^{n}\frac{q^{\frac{n\left(
n-1\right) }{2}}}{\left( \left[ n\right] _{q}\right) !}\frac{v^{n}}{u^{n}}%
t^{n}d_{q}t.  \TCItag{41}
\end{eqnarray}%
\noindent By virtue of $\left( 12\right) ,\left( 41\right) $ yields%
\begin{equation*}
N_{q}\left( D_{q}f\left( t\right) \right) \left( u;v\right) =-f\left(
0\right) +\frac{v}{u}\dint\nolimits_{0}^{\infty }f\left( t\right)
E_{q}\left( -q\frac{v}{u}t\right) d_{q}t.
\end{equation*}%
\noindent Hence,%
\begin{equation*}
N_{q}\left( D_{q}f\left( t\right) \right) \left( u;v\right) =-f\left(
0\right) +\frac{v}{u}N_{q}\left( f\right) \left( u;v\right) .
\end{equation*}%
This completes the proof of the theorem. \\

\noindent Now we extend Theorem 6 to nth derivatives. \\

\noindent \textbf{Theorem 7}. Let $u,v>0$ and $n\in \mathbb{Z}^{+}$. Then, we have%
\begin{equation}
N_{q}\left( D_{q}^{n}f\left( t\right) \right) \left( u;v\right) =\frac{v^{n}%
}{u^{n}}\left( N_{q}\left( f\right) \right) \left( u,v\right)
-\sum_{i=0}^{n-1}\left( \frac{u}{v}\right) ^{n-1-i}D_{q}^{i}f\left( 0\right)
.  \tag{42}
\end{equation}%
\noindent \textbf{Proof}. On account of Theorem 6, we can write%
\begin{eqnarray}
N_{q}\left( D_{q}^{2}f\left( t\right) \right) \left( u;v\right)
&=&D_{q}f\left( 0\right) +\frac{v}{u}N_{q}\left( D_{q}f\right) \left(
u;v\right)  \notag \\
&=&-D_{q}f\left( 0\right) +\frac{v}{u}\left( -f\left( 0\right) +\frac{v}{u}%
\left( N_{q}f\right) \left( u;v\right) \right)  \notag \\
&=&-D_{q}f\left( 0\right) -\frac{v}{u}f\left( 0\right) +\frac{v^{2}}{u^{2}}%
\left( N_{q}f\right) \left( u;v\right) .  \TCItag{43}
\end{eqnarray}%
Proceeding as in $\left( 43\right) $ we obtain%
\begin{equation*}
N_{q}\left( D_{q}^{n}f\left( t\right) \right) \left( u;v\right) =\frac{v^{2}%
}{u^{2}}N_{q}\left( f\right) \left( u;v\right) -\sum_{i=0}^{n-1}\left( \frac{%
v}{u}\right) ^{n-1-i}D_{q}^{i}f\left( 0\right) .
\end{equation*}

\noindent This completes the proof of the theorem.

\section{$N_{q}$ of $q$-Convolutions}

\noindent Let functions $f$ and $g$ be in the form $f\left( t\right)
=t^{\alpha }$ and $g\left( t\right) =t^{\beta -1}$ for $\alpha ,\beta >0$.
We define the $q$-convolution of $f$ and $g$ as%
\begin{equation}
\left( f\ast g\right) _{q}\left( t\right) =\dint\nolimits_{0}^{t}f\left(
\tau \right) g\left( t-q\tau \right) d_{q}t  \tag{44}
\end{equation}%
\noindent \textbf{Theorem 8}. Let $\alpha ,\beta >0$. Then, we have%
\begin{equation}
N_{q}\left( \left( f\ast g\right) _{q}\right) \left( u;v\right)
=u^{2}N_{q}\left( t^{\alpha }\right) \left( u,v\right) N_{q}\left( t^{\beta
-1}\right) \left( u;v\right)  \tag{45}
\end{equation}

\noindent \textbf{Proof} . By aid of $\left( 44\right) $\ and $\left(
14\right) $\ we get%
\begin{equation*}
N_{q}\left( \left( f\ast g\right) _{q}\right) \left( u;v\right) =\frac{%
B_{q}\left( \alpha +1,\beta \right) }{u}\dint\nolimits_{0}^{\infty
}t^{\alpha +\beta }E_{q}\left( -q\frac{vt}{u}\right) d_{q}t.
\end{equation*}%
Hence, by $\left( 28\right) $ we obtain%
\begin{eqnarray}
N_{q}\left( \left( f\ast g\right) _{q}\right) \left( u;v\right) &=&\frac{%
\Gamma _{q}\left( \alpha +1\right) \Gamma _{q}\left( \beta \right) }{\Gamma
_{q}}\frac{u^{\alpha +\beta +1}}{v^{\alpha +\beta +1}}  \notag \\
&=&\Gamma _{q}\left( \alpha +1\right) \Gamma _{q}\left( \beta \right) \frac{%
u^{\alpha +\beta +1}}{v^{\alpha +\beta +1}}.  \TCItag{46}
\end{eqnarray}%
Simple motivation on $\left( 46\right) $ gives%
\begin{equation*}
N_{q}\left( \left( f\ast g\right) _{q}\right) \left( u;v\right) =u^{2}\left(
N_{q}t^{\alpha }\right) \left( u,v\right) \left( N_{q}t^{\beta -1}\right)
\left( u;v\right) .
\end{equation*}%
\noindent The proof is therefore completed. \\

\noindent In similar way, we extend the $q$-convolution to functions of
power series form. \\

\noindent \textbf{Theorem 9}. Let $f\left( t\right) =\sum\limits_{0}^{\infty
}a_{i}t^{\alpha i}$ and $g\left( t\right) =t^{\beta -1}$. Then, we have%
\begin{equation*}
N_{q}\left( \left( f\ast g\right) _{q}\right) \left( u;v\right) =u^{2}\left(
N_{q}f\right) \left( u;v\right) \left( N_{q}g\right) \left( u;v\right) .
\end{equation*}%
\noindent \textbf{Proof}. Under the hypothesis of the theorem and Theorem 8
we write%
\begin{eqnarray*}
N_{q}\left( \left( f\ast g\right) _{q}\right) \left( u;v\right)
&=&\sum\limits_{0}^{\infty }a_{i}N_{q}\left( \left( t^{\alpha i}\ast
t^{\beta -1}\right) _{q}\right) \left( u;v\right)  \\
&=&\sum\limits_{0}^{\infty }a_{i}\left( N_{q}t^{\alpha i}\right) \left(
u;v\right) \left( N_{q}t^{\beta -1}\right) \left( u;v\right)  \\
&=&u^{2}\left( N_{q}f\right) \left( u;v\right) \left( N_{q}g\right) \left(
u;v\right) .
\end{eqnarray*}
\noindent Hence, the proof of the theorem is completed. \\

\section{$N_{q}$ and Heaviside Functions}

\noindent The Heaviside function id defined by%
\begin{equation}
N_{q}\left( \acute{u}\left( t-a\right) \right) =\left\{
\begin{array}{ll}
1 & ,t\geq a \\
0 & ,0\leq t<a%
\end{array}%
\right. ,  \tag{47}
\end{equation}%
where $a$ is a real number. In this part of the paper we merely establish the following
theorem.\\

\noindent \textbf{Theorem 10}. If $\acute{u}$ denotes the heaviside function
and $u,v>0$. Then, we have%
\begin{equation}
N_{q}\left( \acute{u}\left( t-a\right) \right) \left( u;v\right) =\frac{1}{v}%
E_{q}\left( -\frac{v}{u}a\right) .  \tag{48}
\end{equation}%
\noindent \textbf{Proof}. By $\left( 67\right) $ we have%
\begin{equation*}
N_{q}\left( \acute{u}\left( t-a\right) \right) \left( u;v\right) =\frac{1}{u}%
\dint\nolimits_{a}^{\infty }E_{q}\left( -q\frac{v}{u}t\right) d_{q}t.
\end{equation*}%
On account of $\left( 21\right) $ we get%
\begin{eqnarray*}
N_{q}\left( \acute{u}\left( t-a\right) \right) \left( u;v\right)  &=&\frac{1%
}{v}-\frac{1}{u}\dint\nolimits_{0}^{a}\sum\limits_{0}^{\infty }\frac{q^{%
\frac{n\left( n-1\right) }{2}}}{\left( \left[ n\right] _{q}\right) !}\left(
-q\frac{v}{u}t\right) ^{n}d_{q}t \\
&=&\frac{1}{v}-\frac{1}{u}\sum\limits_{0}^{\infty }\left( -1\right) ^{n}%
\frac{q^{\frac{n\left( n-1\right) }{2}}}{\left( \left[ n\right] _{q}\right) !%
}q^{n}\frac{v^{n}}{u^{n}}\dint\nolimits_{0}^{a}t^{n}d_{q}t
\end{eqnarray*}%
Integrating together with simple calculation reveal%
\begin{eqnarray*}
N_{q}\left( \acute{u}\left( t-a\right) \right) \left( u;v\right)  &=&\frac{1%
}{v}-\frac{1}{u}\sum\limits_{0}^{\infty }\frac{q^{\frac{n\left( n-1\right) }{%
2}}}{\left[ n\right] _{q}!}\frac{v^{n}}{u^{n}}\frac{a^{n+1}}{\left[ n+1%
\right] _{q}} \\
&=&\frac{1}{v}+\frac{1}{u}\sum\limits_{0}^{\infty }\left( -1\right) ^{n+1}%
\frac{q^{n+1\frac{n}{2}}}{\left[ n+1\right] _{q}!}\frac{v^{n}}{u^{n}}a^{n+1}
\\
&=&\frac{1}{v}+\frac{1}{v}\sum\limits_{0}^{\infty }\left( -1\right) ^{n+1}%
\frac{q^{\frac{\left( n+1\right) n}{2}}}{\left[ n+1\right] _{q}!}\frac{%
v^{n+1}}{u^{n+1}}a^{n+1}
\end{eqnarray*}%
This can be written as%
\begin{equation*}
N_{q}\left( \acute{u}\left( t-a\right) \right) \left( u;v\right) =\frac{1}{v}%
+\frac{1}{v}\sum\limits_{1}^{\infty }\left( -1\right) ^{m}\frac{q^{m\left(
m-1\right) }}{\left( \left[ m\right] _{q}\right) !}\frac{v^{m}}{u^{m}}a^{m}
\end{equation*}%
Starting the summation from $0$ gives%
\begin{eqnarray*}
N_{q}\left( \acute{u}\left( t-a\right) \right) \left( u;v\right)  &=&\frac{1%
}{v}\sum\limits_{1}^{\infty }\left( -1\right) ^{m}\frac{q^{\frac{m\left(
m-1\right) }{2}}}{\left( \left[ m\right] _{q}\right) !}\frac{v^{m}}{u^{m}}%
a^{m} \\
&=&\frac{1}{v}E_{q}\left( \frac{-v}{u}a\right) .
\end{eqnarray*}

\section{The $q$-Analogue of the $q$-Natural Transform of Second Kind}

\noindent The $q$-analogue of the Natural transform of the second type is
defined over the set $A$,%
\begin{equation*}
A=\left\{ f\left( t\right) \left\vert \exists M,\tau _{1},\tau
_{2}>0,\left\vert f\left( t\right) \right\vert <Me^{t/\tau _{j}},t\in \left(
-1\right) ^{j}\times \left[ 0,\infty \right) \right. \right\}
\end{equation*}%
as%
\begin{equation}
\left( N^{q}f\right) \left( u;v\right) =\frac{1}{u}\dint\nolimits_{0}^{%
\infty }f\left( t\right) e_{q}\left( \frac{-v}{u}t\right) d_{q}t.  \tag{49}
\end{equation}%
The $q$-analogue of the gamma function of the second kind is defined as%
\begin{equation}
\gamma _{q}\left( t\right) =\dint\nolimits_{0}^{\infty }x^{t-1}e_{q}\left(
-x\right) d_{q}x,  \tag{40}
\end{equation}%
and, hence, it follows that%
\begin{equation}
\gamma _{q}\left( 1\right) =1,\gamma _{q}\left( t+1\right) =q^{-t}\left[ t%
\right] _{q}\gamma _{q}\left( t\right) \text{ and }\gamma _{q}\left(
n\right) =q^{\frac{n\left( n-1\right) }{2}}\Gamma _{q}\left( n\right) ,
\tag{41}
\end{equation}%
$\Gamma _{q}$ being the $q$-analogue of gamma function of first kind. \\

\noindent We aim to derive certain results similar to that we have obtained
in the previous sections.\\

\noindent \textbf{Lemma 11}. Let $\alpha >-1$, then we have%
\begin{equation}
\left( N^{q}t^{\alpha }\right) \left( u;v\right) =\frac{u^{\alpha }}{%
v^{\alpha +1}}\gamma _{q}\left( \alpha +1\right) .  \tag{42}
\end{equation}%
\noindent In particular,%
\begin{equation}
\left( N^{q}t^{n}\right) \left( u;v\right) =\frac{u^{n}}{v^{n+1}}q^{\frac{%
-n\left( n-1\right) }{2}}\left( \left[ n_{q}\right] \right) !  \tag{43}
\end{equation}

\noindent \textbf{Proof}. Let $\alpha >-1$, then by change of variables we
have%
\begin{eqnarray*}
\left( N^{q}t^{\alpha }\right) \left( u;v\right) &=&\frac{1}{u}%
\dint\nolimits_{0}^{\infty }t^{\alpha }e_{q}\left( \frac{-vt}{u}\right)
d_{q}t \\
&=&\frac{u^{\alpha }}{v^{\alpha +1}}\dint\nolimits_{0}^{\infty }t^{\alpha
}e_{q}\left( -t\right) d_{q}t
\end{eqnarray*}%
On aid of $\left( 40\right) $, we get%
\begin{equation*}
\left( N^{q}t^{\alpha }\right) \left( u;v\right) =\frac{u^{\alpha }}{%
v^{\alpha +1}}\gamma _{q}\left( \alpha +1\right) .
\end{equation*}%
Proof of the second part of the theorem follows from $\left( 41\right) $. \\

\noindent Hence, we completed the proof of the theorem.

\noindent \textbf{Theorem 12}. Let $a\in \mathbb{R}
, a>0$, then we have%
\begin{equation}
\left( N^{q}e_{q}\left( at\right) \right) \left( u;v\right) =\frac{1}{uv}%
\sum\limits_{0}^{\infty }\frac{a^{n}u^{n}}{v^{n}}q^{\frac{-n\left(
n-1\right) }{2}}.  \tag{44}
\end{equation}

\noindent \textbf{Proof} . By $\left( 13\right) $ we write%
\begin{eqnarray*}
\left( N^{q}e_{q}\left( at\right) \right) \left( u;v\right) &=&\frac{1}{u}%
\dint\nolimits_{0}^{\infty }e_{q}\left( at\right) e_{q}\left( -\frac{v}{u}%
t\right) d_{q}t \\
&=&\frac{1}{u}\sum\limits_{0}^{\infty }\frac{a^{n}}{\left( \left[ n\right]
_{q}\right) !}\dint\nolimits_{0}^{\infty }t^{n}e_{q}\left( -\frac{v}{u}%
t\right) d_{q}t.
\end{eqnarray*}%
By aid of Theorem $11$, the above equation yields%
\begin{equation*}
\left( N^{q}e_{q}\left( at\right) \right) \left( u;v\right) =\frac{1}{uv}%
\sum\limits_{0}^{\infty }\frac{a^{n}u^{n}}{v^{n}}q^{\frac{-n\left(
n-1\right) }{2}}.
\end{equation*}
\noindent This completes the proof of the theorem. \\

\noindent \textbf{Theorem 13}. Let $a>0,a\in
\mathbb{R}
$, then we have%
\begin{equation}
\left( N^{q}E_{q}\left( at\right) \right) \left( u;v\right) =\frac{1}{%
u\left( v-au\right) },\text{ }au<v.  \tag{45}
\end{equation}

\noindent \textbf{Proof}. After some calculations and by using Theorem $11$, we obtain
\begin{eqnarray*}
\left( N^{q}e_{q}\left( at\right) \right) \left( u;v\right) &=&\frac{1}{u}%
\dint\nolimits_{0}^{\infty }E_{q}\left( at\right) e_{q}\left( -\frac{v}{u}%
t\right) d_{q}t \\
&=&\frac{1}{u}\sum\limits_{0}^{\infty }\frac{q^{\frac{n\left( n-1\right) }{2}%
}}{\left( \left[ n\right] _{q}\right) !}a^{n}\dint\nolimits_{0}^{\infty
}t^{n}e_{q}\left( -\frac{v}{u}t\right) d_{q}t \\
&=&\frac{1}{uv}\sum\limits_{0}^{\infty }a^{n}\frac{u^{n}}{v^{n}}.
\end{eqnarray*}%
Since the above series determine a geometric series, we get%
\begin{equation*}
\left( N^{q}e_{q}\left( at\right) \right) \left( u;v\right) =\frac{1}{uv}%
\frac{1}{1-\frac{au}{v}}=\frac{1}{u\left( v-au\right) },au<v.
\end{equation*}%
Hence the theorem is proved.\\

\noindent The $N^{q}$ transform of $\cos ^{q}$ and $\sin ^{q}$ is given as
follows.

\noindent \textbf{Theorem 14}. Let $a>0$, then we have%
\begin{equation}
\left( N^{q}\cos ^{q}\left( at\right) \right) \left( u;v\right) =\frac{1}{u}%
\sum\limits_{0}^{\infty }\left( -1\right) ^{n}a^{2n}\frac{u^{2n}}{v^{2n}}%
q^{-2n\frac{\left( 2n-1\right) }{2}}.  \tag{46}
\end{equation}%
\noindent \textbf{Proof}. Using the definition of $\cos ^{q}$ we write%
\begin{eqnarray*}
\left( N^{q}\cos ^{q}\left( at\right) \right) \left( u;v\right)  &=&\frac{1}{%
u}\dint\nolimits_{0}^{\infty }\cos ^{q}\left( at\right) e_{q}\left( -\frac{v%
}{u}t\right) d_{q}t \\
&=&\frac{1}{u}\sum\limits_{0}^{\infty }\left( -1\right) ^{n}\frac{a^{2n}}{%
\left( \left[ 2n\right] _{q}\right) !}\dint\nolimits_{0}^{\infty
}t^{2n}e^{q}\left( -\frac{v}{u}t\right) d_{q}t \\
&=&\frac{1}{u}\sum\limits_{0}^{\infty }\left( -1\right) ^{n}a^{2n}\frac{%
u^{2n}}{v^{2n}}q^{-2n\frac{\left( 2n-1\right) }{2}}.
\end{eqnarray*}%
This completes the proof of the theorem. \\

\noindent The $N^{q}$ transform of $\sin ^{q}\left( at\right) $ is given as
follows: \\

\noindent \textbf{Theorem 15}. Let $a>0$, then we have%
\begin{equation}
\left( N^{q}\sin ^{q}\left( at\right) \right) \left( u;v\right) =\frac{1}{uv}%
\sum\limits_{0}^{\infty }\left( -1\right) ^{n}a^{2n+1}\frac{u^{2n+1}}{%
v^{2n+1}}q^{-2n\frac{\left( 2n-1\right) }{2}}.  \tag{47}
\end{equation}%
Proof of Theorem 15 follows from similar proof to that of Theorem 14. \\

\section{$N^{q}$ of $q-$Differentiation}

\noindent Before we start investigations, we first assert that
\begin{equation*}
D_{q}e_{q}\left( -\frac{v}{u}t\right) =-\frac{v}{u}e_{q}\left( -\frac{v}{u}%
t\right) .
\end{equation*}%
For further details, we have
\begin{eqnarray*}
D_{q}e_{q}\left( -\frac{v}{u}t\right)  &=&\sum\limits_{0}^{\infty }\frac{%
\left( -1\right) ^{n}}{\left( \left[ n\right] _{q}\right) !}\frac{v^{n}}{%
u^{n}}D_{q}t^{n} \\
&=&\sum\limits_{1}^{\infty }\frac{\left( -1\right) ^{n}}{\left( \left[ n-1%
\right] _{q}\right) !}\frac{v^{n}}{u^{n}}t^{n-1} \\
&=&\sum\limits_{0}^{\infty }\frac{\left( -1\right) ^{n+1}}{\left( \left[ n%
\right] _{q}\right) !}\frac{v^{n+1}}{u^{n+1}}t^{n} \\
&=&-\frac{v}{u}\sum\limits_{0}^{\infty }\frac{\left( -1\right) ^{n}}{\left( %
\left[ n\right] _{q}\right) !}\frac{v^{n}}{u^{n}}t^{n} \\
&=&-\frac{v}{u}e_{q}\left( -\frac{v}{u}t\right) .
\end{eqnarray*}
\noindent This proves the above assertion.\\

\noindent Hence we prove the following theorem. \\

\noindent \textbf{Theorem 16}. Let $u,v>0$, then we have%
\begin{equation}
\left( N^{q}D_{q}f\left( t\right) \right) \left( u;v\right) =-f\left(
0\right) -\dint\nolimits_{0}^{\infty }f\left( qt\right) D_{q}e_{q}\left( -%
\frac{v}{u}t\right) d_{q}t.  \tag{49}
\end{equation}%
\noindent \textbf{Proof.} By $\left( 48\right) $, the above equation gives%
\begin{eqnarray*}
\left( N^{q}D_{q}f\left( t\right) \right) \left( u;v\right)
&=&\dint\nolimits_{0}^{\infty }D_{q}f\left( t\right) e_{q}\left( -\frac{v}{u}%
t\right) d_{q}t \\
&=&-f\left( 0\right) +\frac{v}{u}\dint\nolimits_{0}^{\infty }f\left(
qt\right) e_{q}\left( -\frac{v}{u}t\right) d_{q}t \\
&=&-f\left( 0\right) +\frac{v}{u}q^{-1}\dint\nolimits_{0}^{\infty }f\left(
t\right) e_{q}\left( -\frac{v}{u}q^{-1}t\right) d_{q}t.
\end{eqnarray*}%
By setting variables we have%
\begin{equation*}
\left( N^{q}D_{q}f\left( t\right) \right) \left( u;v\right) =-f\left(
0\right) +\frac{v}{u}q^{-1}\left( N^{q}f\right) \left( q^{-1}v;u\right) .
\end{equation*}%
This completes the proof.\\

\noindent Now, we extend $\left( 49\right) $ to have%
\begin{eqnarray*}
\left( N^{q}D_{q}^{2}f\right) \left( u;v\right)  &=&\left( N^{q}D_{q}\left(
D_{q}f\right) \right) \left( u;v\right)  \\
&=&-D_{q}f\left( 0\right) +\frac{v}{u}q^{-1}\left( N^{q}D_{q}f\right) \left(
q^{-1}v;u\right)  \\
&=&-D_{q}f\left( 0\right) +\frac{v}{u}q^{-1}\left( -f\left( 0\right) +\frac{v%
}{u}q^{-1}\left( N^{q}f\right) \left( q^{-2}v;u\right) \right)  \\
&=&-D_{q}f\left( 0\right) -\left( \frac{v}{u}\right) q^{-1}f\left( 0\right)
+\left( \frac{v}{u}\right) q^{2-2}\left( N^{q}f\right) \left(
q^{-2}v;u\right) .
\end{eqnarray*}%
Proceeding to nth derivatives, we get%
\begin{equation*}
\left( N^{q}D_{q}^{n}f\right) \left( u;v\right) =\left( \frac{v}{u}\right)
^{n}q^{-n}\left( N^{q}f\right) \left( q^{-n}v;u\right)
-\sum\limits_{i=0}^{n-1}\left( \frac{v}{u}\right) ^{n-1-i}D_{q}^{i}f\left(
0\right) .
\end{equation*}%
This completes the proof of the theorem.

\section*{Competing interests } The authors declare that they have no competing
interests.

\section*{Author's contributions}
    All the authors jointly worked on deriving the results and approved the final manuscript.


\begin{thebibliography}{99}

\bibitem{2} W. H. Abdi, On q-Laplace transforms, Proc. Nat. Acad. Sci. India
29 (1961), 389--408

\bibitem{5} D. Albayrak, S. D. Purohit , Faruk Ucar, On q-Sumudu transforms
of certain $q$-Polynomials, Filomat 27(2)(2013), 413--429.

\bibitem{6} D. Albayrak, S. D. Purohit and F. Ucar, On $q$-analogues of
Sumudu transform, An. St. Univ. Ovidius Constant, 21(1)(2013), 239--260

\bibitem{17} S. K. Q. Al-Omari, On the application of the Natural
transforms, Inter. J. Pure Appl. Math. 85(4), (2013), 729--744.

\bibitem{14} F. B. M. Belgacem and R. Silambarasan, Theoretical
investigations of the Natural transform, Progress In Electromagnetics
Research Symposium Proceedings, Suzhou, China, Sept. (2011), 12--16.

\bibitem{19} F. B. M. Belgacem and R. Silambarasan, Advances in the Natural
transform, AIP Conf. Proc. 1493, 106 (2012), doi: 10.1063/1.4765477.

\bibitem{15} F. B. M. Belgacem and R. Silambarasan, Maxwell's equations
solutions through the Natural transform,\ Mathematics in Engineering,
Science and Aerospace 3(3), (2012), 313--323.

\bibitem{8} A. Fitouhi and N. Bettaibi. Wavelet transforms in quantum
calculus. Journal of Nonlinear Mathematical Physics 13(3), (2006), 492--506

\bibitem{9} A. Fitouhi, N\'{e}ji Bettaibi. Applications of the Mellin
transform in quantum calculus, J. Math. Anal. Appl. 328 (2007), 518--534.

\bibitem{10} G. Gasper, M. Rahmen, Basic hypergeometric series, Encyclopedia
Math. Appl., 35, Cambridge Univ. Press,Cambridge, UK, 1990.

\bibitem{18} W. Hahn, Beitrage zur theorie der heineschen reihen, die 24
Integrale der hypergeometrischen $q$-Diferenzengleichung, das $q$-analog on
der Laplace transformation, Math. Nachr. 2 (1949), 340--379.

\bibitem{11} F. H. Jackson, On a q-Definite integrals, Quart. J. Pure and
Appl. Math. 41 (1910), 193--203.

\bibitem{12} V. G. Kac, P. Cheung, Quantum calculus, Universitext,
Springer-Verlag, New York, 2002

\bibitem{13} Z. H. Khan and W. A. Khan, N-transform properties and
applications,\ NUST Jour. of Engg. Sciences, (1)(2008), 127--133.

\bibitem{3} S. D. Purohit and S. L. Kalla, On $q$-Laplace transforms of the $%
q$-Bessel functions, Fract. Calc. Appl. Anal. 10(2)(2007), 189--196.

\bibitem{16} R. Silambarasan and F. B. M. Belgacem, Applications of the
Natural transform to Maxwell's equations, Progress In Electromagnetics
Research Symposium Proceedings, Suzhou, China, Sept. (2011),
12--16.

\bibitem{4} F. Ucar and D. Albayrak, On $q$-Laplace type integral operators
and their applications, Journal of Difference Equations and Applications,
iFirst article, (2011), 1--14.

\bibitem{7} R. Yadav and S. D. Purohit, On Applications of Weyl fractional $q
$-integral operator to generalized basic hypergeometric functions, Kyungpook
Math. J. 46, (2006), 235--245.
\end{thebibliography}
\end{document}